\documentclass[12pt,a4paper,reqno,draft]{amsart}
\usepackage{amsmath,amsfonts,amsthm,amssymb,color}
\usepackage{cmmib57}
\usepackage{exscale}
\usepackage{amscd}
\usepackage{latexsym}
\usepackage{graphicx}
\usepackage[T1]{fontenc}
\usepackage[latin1]{inputenc}
\usepackage{pdfsync}

\newtheorem{thm}{Theorem}[section]
\newtheorem{cor}[thm]{Corollary}
\newtheorem{lem}[thm]{Lemma}
\newtheorem{prop}[thm]{Proposition}
\newtheorem{defn}{Definition}[section]

\theoremstyle{remark}

\begin{document}

\title[Inviscid limit for 2D SNSE]{Inviscid limit for 2D stochastic \\ Navier-Stokes equations}

\author{Fernanda Cipriano \and Iván Torrecilla}

\address{{\it Fernanda Cipriano:} {\rm Grupo de Física Matemática e Dep. de Matemática FCT, Universidade Nova de Lisboa, Av. Prof. Gama Pinto, 2, PT-1649-003 Lisboa, Portugal}. \newline {\it Email: }{\tt cipriano@cii.fc.ul.pt}
\vspace{0.5cm}
\newline $\mbox{ }$\hspace{0.1cm} {\it Iván Torrecilla:} {\rm Grupo de Física Matematica, Instituto para a Investigaç\~{a}o Interdisciplinar da Universidade de Lisboa, Av. Prof. Gama Pinto, 2, PT-1649-003 Lisboa, Portugal}. \newline {\it Email: }{\tt itorrecillatarantino@gmail.com}.}

\begin{abstract}
We consider stochastic Navier-Stokes equations in a 2D-bounded domain with
 the Navier with friction boundary condition.
 We establish the existence and the uniqueness of the solutions  and study the vanishing viscosity limit. More precisely, we prove that
solutions of stochastic Navier-Stokes equations converge, as the viscosity goes to zero, to solutions of the corresponding stochastic Euler
equations.
\end{abstract}


\date{\today}
\maketitle


\medskip


\section{Introduction}
\label{sec0}
The study of the inviscid limit of the solutions of the Navier-Stokes equations is a classical issue in fluid mechanics.
The knowledge of the behavior of the solutions for small viscosities (very high Reynolds numbers) is crucial to
understand the turbulence phenomena. The mathematical resolution of the
inviscid limit problem should have strong consequences in many branches of
engineering (technology involving heat and mass transfer),  as aircraft
production, turbine blades, nanotechnology, etc..

\smallskip

The investigation of this problem for domains without boundary was performed, for instance, in \cite{B10}, \cite{CF88},  \cite{CW95}, \cite{K72}. When the Navier-Stokes equations are considered with a stochastic random force,
the  inviscid limit of its solutions is studied in \cite{BP01}.

\smallskip

In the case of bounded domains, the Navier-Stokes equations should be supplemented with boundary conditions. The most studied and  widely accepted is the Dirichlet boundary condition which prescribes the value of the velocity field on the surface boundary. In the presence of the impermeable boundary, the normal  and the tangential components of the velocity are assumed to be zero on the surface. For the Euler equations it is just required that the velocity field be tangent to the boundary. In the vanishing viscosity strong boundary layers arise, which are very difficult to treat and
the inviscid limit remains an open problem.
 Just
partial results have been obtained (see \cite{SC98}, \cite{TW02}).
Other physical meaningful boundary condition is  the so called Navier slip boundary condition. this boundary condition was initially introduced  by Navier \cite{N1827} in 1827  and due to recent experimental results (see \cite{BFN10}, \cite{JM01}, \cite{PT06}),  has been renewed interest in this boundary condition.

\smallskip

To be more precise,
we suppose that $\mathcal{O}$ is a bounded simply connected domain in $\mathbb{R}^2$ with  boundary  $\Gamma$ sufficiently regular. The Navier slip with friction boundary condition, for the Navier-Stokes equations,  is written by
\begin{equation}
\label{eqNSBC}
2D(u)\mathbf{n}\cdot \mathbf{t}+\alpha u\cdot \mathbf{t}=0
 \quad\mbox{ on }\quad ]0,T[\times\Gamma
\end{equation}
where $\displaystyle{
D(u)=\frac{1}{2}
\bigl(\nabla u+(\nabla u)^T \bigr)}$ is the rate-of-strain tensor;
  $\mathbf{n}$ and $\mathbf{t}$ are the unit exterior normal and the unit tangent vector, respectively, to $\Gamma$, $\{\mathbf{n},\mathbf{t}\}$ being a direct basis; and $"\cdot"$ defines the scalar product on $\mathbb{R}^2$.
 Here the tangent component of  the fluid velocity at the boundary, rather  than being fixed, is proportional to the tangential stress.
The normal component of the  fluid velocity at the boundary is zero and corresponds to the impermeability of the boundary:
\begin{align}
u \cdot \mathbf{n} =0 \quad\mbox{ on }\quad ]0,T[ \times \Gamma.
\label{eqFBC}
\end{align}

\smallskip

The key feature of this  boundary condition \eqref{eqNSBC}-\eqref{eqFBC} is that it can be expressed in terms of the vorticity $\xi$ of the vector field $u$ as
\begin{equation}
\label{eqFBC1}
\xi(u)=(2\kappa-\alpha)u \cdot\mathbf{t}  \quad\text{ and }\quad \quad u \cdot \mathbf{n} =0 \quad\text{ on } ]0,T[\times\Gamma,
\end{equation}
which permits  to handle the vorticity formulation of the Navier-Stokes equations.
The coefficient $\alpha  $ is a  known function describing physical properties and $\kappa$ is the curvature of the boundary.

\smallskip

A particular case of this boundary condition, with  $\alpha=2k,$
\begin{equation}
\xi =0\quad \quad \text{on }\quad \quad ]0,T[\times\Gamma \label{HBCw-nsw}
\end{equation}%
was considered in \cite{L96}; where an energy type estimate for $%
\xi $ was established, allowing to prove the
 convergence of the solutions of the Navier-Stokes equations to solutions of the Euler equations.
 This boundary condition is also known as the Lions boundary condition or free boundary condition.
  Besides its mathematical importance this particular boundary condition do not permits the creation of the vorticity on the boundary.
 The deterministic methods were extended in \cite{BF99} and \cite{BP01}, to obtain some well posedness results  for 2D stochastic Euler equations.
 In both articles, the stochastic Euler equations are regularized by the corresponding viscous stochastic Navier-Stokes equations supplemented with the Lions boundary condition; the zero-viscosity limit provides the solution
 for the stochastic Euler equations.
  In \cite{BF99} is considered an additive noise and the inviscid limit is a strong solution (in the probability sense) of stochastic Euler equations. Moreover, a uniqueness result is established if the initial vorticity belongs to $L^{\infty }.$
  In \cite{BP01}, a less regular multiplicative noise is considered  and the inviscid limit gives a martingale solution to the stochastic Euler equations.
  More recently,  \cite{BM10} handled this particular case of the Navier slip boundary condition for stochastic Navier-Stokes equations with a multiplicative noise and studied the viscous limit  using the large deviations techniques, taking the square root of the viscosity in front of the noise.

\smallskip

 In the deterministic framework, the study of the inviscid limit for the solutions of the Navier-Stokes equations
 with the physical Navier slip boundary conditions \eqref{eqNSBC} has been greatly developed.
 In  \cite{CMR98},  the solvability of the Navier Stokes equations
with the boundary condition \eqref{eqNSBC} was established in the class of $L^{\infty }$-bounded
vorticity. It was  also proved that the vanishing viscosity limit is
well described by Euler equations.
Later on, in \cite{LFLP05}, this result was generalized for the class of
$L^{p}$-bounded vorticity with $p>2$.
 A rate of the vanishing viscous
convergence of solutions of the Navier-Stokes equations to solutions
 of the Euler equations, in the class of almost $L^{\infty }$-bounded vorticity,  was obtained
in \cite{Ke06}.

\smallskip

 In the present  work we consider stochastic Navier-Stokes equations,  with an additive noise, on a bounded domain of $\mathbb{R}^2$, subjected to the Navier slip with friction boundary condition \eqref{eqNSBC}-\eqref{eqFBC}, which provides creation of the vorticity on the boundary  proportional to the tangential velocity, and  tackle the problem of the inviscid limit.
In some sense, our result is the probability  counterpart of the deterministic result obtained in
 \cite{LFLP05}.

\smallskip

The article is organized as follows:\ \ in the section 2, we introduce the functional spaces,
 construct the appropriate Wiener process and  state the main result Theorem \ref{main-result}. In the Section 3, we deduce  the $L^2$ a priori estimates for the viscous solutions independent of $\nu.$   The section 4 contains
  the relevant  $L^p$ a priori estimates for the viscous vorticity  independent  of $\nu$. These estimates
  permit to establish the well posedness for the Navier-Stokes equations with the Navier boundary condition.
  In the last section we obtain crucial
  path-wise estimates independent of the viscosity, that
    allow to  establish the inviscid limit.

\smallskip

\section{Velocity equations with additive noise}
\label{sec1}
We consider the following stochastic Navier-Stokes equations in dimension 2:
\begin{equation}
\label{NS:eq}
\left\{\renewcommand{\arraystretch}{1.8}
\begin{array}{ll}
\frac{\partial u^\nu(t)}{\partial t}-\nu \Delta u^\nu(t)+(u^\nu(t) \cdot\nabla )u^\nu(t)+\nabla p(t)= f(t)+\sqrt{\mathcal{Q}}\,\dot{W} &\text{ in }]0,T[\times \mathcal{O},    \\
{\rm div\,} u^\nu=0 &\text{ in }]0,T[\times \mathcal{O}, \\
u^\nu(0)=u_0 &\text{ in } \mathcal{O},\\
u^\nu\cdot \mathbf{n} =0 &\text{ on }]0,T[\times \Gamma,\\
2D(u^\nu)\mathbf{n}\cdot \mathbf{t}+\alpha u^\nu\cdot \mathbf{t}=0 &\text{ on }]0,T[\times \Gamma
\end{array}
\right.
\end{equation}
where $\nu>0$ is the coefficient of kinematic viscosity, $\Delta$ denotes the Laplacian, $\nabla$ denotes the gradient, ${\rm div}\, u^\nu = \nabla\cdot u^\nu= \sum_{i=1}^2 \partial_i u^{\nu,i}$, $\alpha(x)$ is a given positive twice continuously differentiable function defined on $\Gamma\doteq \partial \mathcal{O}$, $u^\nu$ is the velocity and $p$ is the pressure. $f(t,x)$ is a given deterministic force and $\sqrt{\mathcal{Q}}\,\dot{W}$ is the formal derivative of a Gaussian random field  in time and correlated in space that will be set below.

\smallskip

We introduce the following Hilbert spaces
\begin{align*}
H&=\left\{v\in \left[L^2(\mathcal{O})\right]^2: \nabla\cdot v=0 \text{ in }\mathcal{O}\text{ and }v\cdot \mathbf{n}=0 \text{ on }\Gamma\right\},\\
V&=\left\{v\in \left[H^1(\mathcal{O})\right]^2: \nabla\cdot v=0 \text{ in }\mathcal{O}\text{ and }v\cdot \mathbf{n}=0 \text{ on }\Gamma\right\},\\
\mathcal{W}&=\left\{v\in V\cap\left[H^2(\mathcal{O})\right]^2: 2D(u^\nu)
\mathbf{n}\cdot \mathbf{t}+\alpha u^\nu\cdot \mathbf{t}=0 \text{ on }\Gamma\right\}.
\end{align*}
It can be verified (see Lemma 2.1 of \cite{CMR98}) that
\begin{align*}
\mathcal{W}&=\left\{v\in V\cap\left[H^2(\mathcal{O})\right]^2: {\rm curl}\,v=(2\kappa-\alpha)v\cdot \mathbf{t} \text{ on }\Gamma\right\}
\end{align*}
where $\kappa$ denotes the curvature of $\Gamma$. We recall that ${\rm curl}\,v=\partial_1 v^2-\partial_2 v^1$.

\smallskip

We consider on $H$ the
$L^2$- inner product and norm that we denote by $\left\langle \cdot,\cdot\right \rangle $ and  $\|\cdot\|_{L^2}$.  $V$ is endowed with
the inner product $$\left\langle u,v\right \rangle_V=\left\langle \nabla u,
\nabla v\right \rangle$$ and the associated norm $\|\cdot\|_V$. We recall that, from the Poincar\'e's inequality, this norm is equivalent to the $H^1$-norm

\smallskip

Let us denote by $V^\prime$  the topological dual of $V$ and by
$\left \langle \cdot ,\cdot
\right\rangle_{V^{\prime}, V}$
 the corresponding  duality. We define the operator
$\mathcal{A}:V\rightarrow V^\prime$ by
\begin{equation}
\label{A:def}
\left \langle \mathcal{A}u,v\right \rangle_{V^\prime, V}=\int_\mathcal{O} \nabla u\cdot \nabla v-\int_\Gamma (\kappa-\alpha)u\cdot v,
\end{equation}
for all $u,v\in V$.
Since
$$
\bigl|\left \langle \mathcal{A}u,v\right \rangle_{V^\prime, V}\bigr|\leq C\|u\|_V\|v\|_V
$$
 $\mathcal{A}$ is a continuous operator form $V$ to $V^\prime.$
 Moreover $\mathcal{A}: \mathcal{W}\to H$ coincides  with the stokes operator $-P_H \Delta $, where $P_H$ denotes the Leray projector. More precisely we have
$$\left \langle \mathcal{A}u,v\right \rangle_{V^\prime,V}=\left\langle -\Delta u,v\right \rangle, \quad u\in \mathcal{W},
v\in V.
$$

\smallskip

We also define $\mathcal{B}:V\rightarrow V^\prime$ as $\mathcal{B} (u)=(u \cdot\nabla )u$, that is,
\begin{equation}
\label{B:def}
\left \langle \mathcal{B}(u),v\right \rangle=\int_\mathcal{O} (u \cdot\nabla )u\cdot  v,
\end{equation}
for all $u,v\in V$.
\smallskip

\smallskip
From Lemma 2.2 of \cite{CMR98}, there exists a basis $\{v_k\} \subset \mathcal{W}$ for $V$, of eigenfunctions of the operator $\mathcal{A}$, being simultaneously  an orthonormal basis for   $H.$ The corresponding sequence $\{\lambda_k\}$ of  eigenvalues verifies $\lambda_k>0$, $\forall k\in \mathbb{N}$ and $\lambda_k \to \infty$
as $k\to \infty.$  Henceforth we shall consider this basis.

To be more specific, we shall take in the following $\mathcal{Q}=\mathcal{A}^{-2m}$, where $m\in \mathbb{N}$ will be fixed later and $W(t)=\sum_{k=1}^\infty \beta_k(t)v_k$, $t\geq 0$. Here $\{\beta_k\}$ denotes a sequence of standard Brownian motion mutually independent defined on a filtered probability space  $(\Omega,\mathcal{F},\mathbf{P},\{\mathcal{F}_t\}_{t\geq 0})$. In fact,
$$
\sqrt{\mathcal{Q}}\,W(t)=\sum_{k=1}^\infty \beta_k(t)\sqrt{\mathcal{Q}}\,v_k=\sum_{k=1}^\infty \lambda_k^{-m}v_k\beta_k(t)
$$
is a $H$-valued centered Wiener process on $(\Omega,\mathcal{F},\mathbf{P})$, with covariance $\mathcal{Q}$ in $H$. We take $m\in\mathbb{N}$ such that
\begin{equation}
\label{Q:cond}
\mathcal{M}:=\sum_{k=1}^\infty \lambda_k^{-2m+3}<\infty.
\end{equation}
Then, with this choice of $m$ we have that $\mathcal{Q}$ is an operator of trace class. We denote the trace of $\mathcal{Q}$ by ${\rm tr}(\mathcal{Q})\doteq \sum_{k=1}^\infty \langle \mathcal{Q}v_k,v_k\rangle=\sum_{k=1}^\infty \lambda_k^{-2m}$. Let us mention  that a similar noise was considered  in \cite{BDP08}.

In terms of $\mathcal{A}$, $\mathcal{B}$ and $f$ we can write Equation (\ref{NS:eq}) as the following stochastic evolution equation in $V'$:
\begin{equation}
\label{NS:eqevol}
\left\{\renewcommand{\arraystretch}{1.8}
\begin{array}{ll}
du^\nu=F(t,u^\nu(t))\,dt+\sqrt{\mathcal{Q}}\,dW(t) &\text{ in }]0,T[\times \mathcal{O},    \\
u^\nu(0)=u_0 &\text{ in } \mathcal{O},
\end{array}
\right.
\end{equation}
where $F(t,u^\nu)=f-\nu\mathcal{A}u^\nu-\mathcal{B}(u^\nu)$.

\begin{defn}
\label{definicao}
Given $u_0\in L^2\left(\Omega;H\right)$, an adapted stochastic process $u^\nu$ with sample paths in
$\mathcal{C}([0,T]; H)\cap L^2(0,T;V)$ is said a weak solution of the stochastic Navier-Stokes equation (\ref{NS:eqevol}) if
\begin{equation}
\label{NS:eqweak2}
\left\langle u^\nu(t),v\right\rangle=
\left\langle u_0,v\right\rangle
+
\int_0^t\left\langle F(s,u^\nu(s)),v\right\rangle\,ds+\int_0^t\left\langle \sqrt{\mathcal{Q}}\,dW(s),v\right \rangle,
\end{equation}
in $]0,T[$, for all $v\in V$ and a.e.
 $\omega\in \Omega$.

\end{defn}

\smallskip

For  the viscosity equal to zero we consider the stochastic two-dimensional Euler equations
\begin{equation}
\label{E:eq}
\left\{\renewcommand{\arraystretch}{1.8}
\begin{array}{ll}
\frac{\partial u(t)}{\partial t}+(u(t) \cdot\nabla )u(t)+\nabla p(t)= f(t)+\sqrt{\mathcal{Q}}\,\dot{W} &\text{ in }]0,T[\times \mathcal{O},    \\
{\rm div\,} u=0 &\text{ in }]0,T[\times \mathcal{O}, \\
u(0)=u_0 &\text{ in } \mathcal{O},\\
u\cdot \mathbf{n} =0 &\text{ on }]0,T[\times \Gamma
\end{array}
\right.
\end{equation}
which can be written in terms of the operator $\mathcal{B}$ and $f$  by the following stochastic evolution equation in $V'$:
\begin{equation}
\label{E:eqevol}
\left\{\renewcommand{\arraystretch}{1.8}
\begin{array}{ll}
du(t)=\left\{f(t)-\mathcal{B}(t,u(t))\right\}\,dt
+\sqrt{\mathcal{Q}}\,dW(t) &\text{ in }]0,T[\times \mathcal{O},    \\
u^\nu(0)=u_0 &\text{ in } \mathcal{O}.
\end{array}
\right.
\end{equation}

\smallskip

The main result of this article is the following:
\begin{thm}
\label{main-result}Let $T>0$, $\nu_0>0$ and $p>2$. Suppose that $f\in L^2\left(0,T;H\right)$, ${\rm curl}\,f\in L^1\left(0,T;L^p(\mathcal{O})\right)$, $u_0\in L^p(\Omega;H)$ and ${\rm curl}\,u_0\in L^p(\Omega;L^p(\mathcal{O}))$. Then we have:
\begin{itemize}
\item[(i)]For any $\nu\in]0,\nu_0]$, there exists a unique  weak solution $u^\nu $ of the stochastic Navier-Stokes equation (\ref{NS:eqevol})
 such that
\begin{align*}
&u^\nu\in L^p\left(\Omega;\mathcal{C}([0,T];H)\right)\cap L^2\left(\Omega;L^2(0,T;V)\right)\cap [L^4\left(]0,T[\times \mathcal{O}\times \Omega\right)]^2,
\\
&{\rm curl}\,u^\nu\in L^2(\Omega;
L^{\infty}(0,T; L^p(\mathcal{O})).
\end{align*}
\smallskip

\item[(ii)] In addition, if
  ${\rm curl}\,f\in L^1\left(0,T;L^\infty(\mathcal{O})\right)$,
 there exists  a measurable stochastic process $u$ that is a  solution of the incompressible 2D stochastic
 Euler equation \eqref{E:eqevol},  in the  sense that
\begin{align}
\label{E:eqweak:time}
\left\langle u(t),v\right\rangle &= \left\langle u_0,v
\right\rangle
-\int_0^t\left\langle \mathcal{B}(u(s)),v\right\rangle\,ds+\int_0^t\left\langle f(s),v\right\rangle\,ds\notag\\
&\quad+\int_0^t\left\langle \sqrt{\mathcal{Q}}\,dW(s),v\right \rangle
\end{align}
for all $v\in V$ and $\mathbb{P}$-a.e.  $\omega\in \Omega.$
 Furthermore, taking ${\rm curl}\,u_0\in L^p(\Omega;L^\infty(\mathcal{O}))$,
for $\mathbb{P}$-a.e. $\omega\in \Omega$
$$
u^\nu(\omega)\to u(\omega) \quad \text{strongly in  }\;\;
\mathcal{C}\left([0,T];H\right), \text{ as } \nu\to 0.
$$
\end{itemize}
\end{thm}

\smallskip

\section{$L^2$ a priori estimates for the velocity and  solvability of the Navier-stokes equations}
\label{sec3}
We consider the following Faedo-Galerkin approximations of Equation (\ref{NS:eqevol}). Let $H_n\doteq {\rm span}\,\{v_1,\ldots,v_n\}$  and define $u^\nu_n$ as the solution of the following stochastic differential equation:

\smallskip

For each $v\in H_n$,
\begin{equation}
\label{NS:sde}
d\langle u^\nu_n(t),v\rangle=\langle F(t,u^\nu_n(t)),v \rangle\,dt+\langle \sqrt{\mathcal{Q}}\,dW(t),v\rangle,
\end{equation}
with $u^\nu_n(0)=\sum_{k=1}^n \langle u_0,v_k\rangle v_k$.

\smallskip

Notice that Equation (\ref{NS:sde}) defines a system of stochastic ordinary differential equations in $\mathbb{R}^n$ with locally Lipschitz coefficients. Therefore, we need some a priori estimate to prove the global existence of a solution $u^\nu_n(t)$ as an adapted process in the space $\mathcal{C}([0,T];H_n)$.

\smallskip

\begin{prop}
\label{energy:est}
Let $T>0$ and $\nu_0>0$. Suppose that $f\in L^1(0,T;H)$ and $u_0\in L^2(\Omega;H)$. Let $u^\nu_n(t)$ be an adapted process in the space $\mathcal{C}([0,T];H_n)$ solution of Equation (\ref{NS:sde}). Then
\begin{align}
\label{apriori-est}
&\sup_{0<\nu\leq \nu_0}\sup_{n}\left\{
\mathbb{E}\left(\sup_{0\leq r \leq T}\|u^\nu_n(r)\|_{L^2}^2\right)+\nu\int_0^T\mathbb{E}\left(\|u^\nu_n(s)\|_V^2\right)\,ds\right\}\notag\\
&\leq C(f,\mathcal{Q},\nu_0)\left(\mathbb{E}\left(\|u_0\|_{L^2}^2\right)+1\right).
\end{align}
Furthermore  we have
\begin{align}
\label{energy:eqlim}
&\|u^\nu_n(t)\|_{L^2}^2+2\nu\int_0^{t}\|\nabla u^\nu_n(s)\|_{L^2}^2\,ds\notag\\
&=\|u^\nu_n(0)\|_{L^2}^2
+2\nu\int_0^{t}\left(\int_{\Gamma}(\kappa-\alpha)u^\nu_n(s)\cdot u^\nu_n(s)\,d\mathcal{S}\right)ds+
2\int_0^{t}\langle  f(s),u^\nu_n(s)\rangle \,ds\notag\\
&\quad +2\int_0^{t}\langle \sqrt{\mathcal{Q}}\,dW(s),u^\nu_n(s)\rangle +\int_0^{t} {\rm tr}(\mathcal{Q})\,ds,
\end{align}
\end{prop}

\smallskip

\textbf{Proof}.
{ For each $N\in \mathbb{N}$, let us consider the stopping time
$\tau_N=\inf\{t\geq 0:\|u_n^\nu(t)\|_{L^2}\geq N\}\wedge T$}.
From It\^o's formula 
\begin{align}
\label{NS:ito1}
\|u^\nu_n({t\wedge\tau_N})\|_{L^2}^2&=\|u^\nu_n(0)\|_{L^2}^2+2\int_0^{t\wedge\tau_N}\langle   F(s,u^\nu_n(s)),u^\nu_n(s)\rangle \,ds\notag\\
&\quad +2\int_0^{t\wedge\tau_N}\langle \sqrt{\mathcal{Q}}\,dW(s),u^\nu_n(s)\rangle +\int_0^{t\wedge\tau_N} {\rm tr}(\mathcal{Q})\,ds.
\end{align}

\smallskip

Applying to (\ref{NS:ito1}) the definition of operator $\mathcal{A}$ (\ref{A:def}) and of operator $\mathcal{B}$ (\ref{B:def}), respectively, and integration by parts formula, we obtain expression (\ref{energy:eq}).

{\begin{align}
\label{energy:eq}
&\|u^\nu_n(t\wedge\tau_N)\|_{L^2}^2+2\nu\int_0^{t\wedge\tau_N}\|\nabla u^\nu_n(s)\|_{L^2}^2\,ds\notag\\
&=\|u^\nu_n(0)\|_{L^2}^2
+2\nu\int_0^{t\wedge\tau_N}\left(\int_{\Gamma}(\kappa-\alpha)u^\nu_n(s)\cdot u^\nu_n(s)\,d\mathcal{S}\right)ds+
2\int_0^{t\wedge\tau_N}\langle  f(s),u^\nu_n(s)\rangle \,ds\notag\\
&\quad +2\int_0^{t\wedge\tau_N}\langle \sqrt{\mathcal{Q}}\,dW(s),u^\nu_n(s)\rangle +\int_0^{t\wedge\tau_N} {\rm tr}(\mathcal{Q})\,ds,
\end{align}

}

\smallskip

Moreover,
\begin{align}
\label{NS:v1}
&\int_{\Gamma}(\kappa-\alpha)u^\nu_n(s)\cdot u^\nu_n(s)\,d\mathcal{S}\leq \sup_\Gamma |\kappa-\alpha|\|u^\nu_n(s)\|_{L^2(\Gamma)}^2\notag\\
&\leq \sup_\Gamma |\kappa-\alpha| C(\mathcal{O})\frac{1}{\sqrt{2\varepsilon}}\|u^\nu_n(s)\|_{L^2}\sqrt{2\varepsilon}\|\nabla u^\nu_n(s)\|_{L^2}\notag\\
&\leq \varepsilon \|\nabla u^\nu_n(s)\|_{L^2}^2+C(\varepsilon)\|u^\nu_n(s)\|_{L^2}^2,
\end{align}
where $C(\varepsilon)\doteq \sup_\Gamma |\kappa-\alpha|^2 C(\mathcal{O})^2\frac{1}{4\varepsilon}$.

\smallskip

The application of Cauchy-Schwarz's inequality gives
\begin{align}
\label{NS:v2}
|\langle f(s),u^\nu_n(s)\rangle|\leq \|f(s)\|_{L^2}\|u^\nu_n(s)\|_{L^2}\leq \|f(s)\|_{L^2}\left(1+\|u^\nu_n(s)\|_{L^2}^2\right).
\end{align}

\smallskip

Applying Burkholder-Davis-Gundy's inequality
\begin{align}
\label{NS:v3}
\mathbb{E}\left|\sup_{0\leq r\leq t} \left\{2\int_0^{r\wedge\tau_N}\langle \sqrt{\mathcal{Q}}\,dW(s),u^\nu_n(s)\rangle\right\}\right|&\leq 2C_1{\rm tr}(\mathcal{Q})^{1/2}\mathbb{E}\left(\int_0^{t\wedge\tau_N}\|u^\nu_n(s)\|_{L^2}^2\,ds\right)^{1/2}\notag\\
&\leq C_1 {\rm tr}(\mathcal{Q}) + C_1\mathbb{E}\left(\int_0^{t\wedge\tau_N}\|u^\nu_n(s)\|_{L^2}^2\,ds\right).
\end{align}

\smallskip

Using expression (\ref{energy:eq}), estimates (\ref{NS:v1}), (\ref{NS:v2}) and (\ref{NS:v3}), and $\nu\leq\nu_0$, we obtain
\begin{align}
\label{NS:gronbell0}
&\mathbb{E}\left(\sup_{0\leq r \leq t}\|u^\nu_n({r\wedge\tau_N})\|_{L^2}^2\right)+2\nu(1-\varepsilon)
\int_0^{t\wedge\tau_N}\mathbb{E}\left(\|u^\nu_n(s)\|_V^2\right)\,ds\notag\\
&\leq \mathbb{E}\left(\|u_0\|_{L^2}^2\right)+2\nu_0 (C(\varepsilon)+1)\int_0^{t\wedge\tau_N}\mathbb{E}
\left(\sup_{0\leq r\leq s}\|u^\nu_n({r\wedge\tau_N})\|_{L^2}^2\right)\,ds\notag\\
&\quad+\int_0^{t\wedge\tau_N}\left[{\rm tr}(\mathcal{Q})+2\|f(s)\|_{L^2}\right]\,ds+2\int_0^{t\wedge\tau_N}\|f(s)\|_{L^2}\mathbb{E}\left(\sup_{0\leq r \leq s}\|u^\nu_n({r\wedge\tau_N})\|_{L^2}^2\right)\,ds\notag\\
&\quad+C_1 {\rm tr}(\mathcal{Q}) + C_1\int_0^{t\wedge\tau_N}\mathbb{E}\left(\sup_{0\leq r \leq s}\|u^\nu_n({r\wedge\tau_N})\|_{L^2}^2\right)\,ds.
\end{align}

Finally, in (\ref{NS:gronbell0}) set $\varepsilon=1/2$ and apply Gronwall-Bellman inequality (see pp. 651-652 in \cite{Kh02}) to
\begin{align*}
X(t)&\doteq \mathbb{E}\left(\sup_{0\leq r \leq t}\|u^\nu_n({r\wedge\tau_N})\|_{L^2}^2\right),\quad Y(t)\doteq \nu\int_0^{t\wedge\tau_N}\mathbb{E}\left(\|u^\nu_n(s)\|_V^2\right)\,ds,\\
Z(t)&\doteq \mathbb{E}\left(\|u_0\|_{L^2}^2\right)+C_1{\rm tr}(\mathcal{Q})+\int_0^{t\wedge\tau_N} \left[{\rm tr}(\mathcal{Q})+2\|f(s)\|_{L^2}\right]\,ds,\\
\lambda(t)&\doteq Z(t)-Y(t),\\
\varphi(t)&\doteq 2\nu_0(C(1/2)+1)+C_1+2\|f(s)\|_{L^2}\geq 0 \text{ with }\mathcal{K}\doteq \int_0^T\varphi(s)\,ds.
\end{align*}
Then,
\begin{align}
\label{eq:ee}
X(t)+Y(t)\leq Z(T)(1+\mathcal{K}e^{\mathcal{K}})\leq C(f,\mathcal{Q},\nu_0)\left(\mathbb{E}\left(\|u_0\|_{L^2}^2\right)+1\right),
\end{align}
uniformly in ${N}$, $n$, and $\nu\leq \nu_0$.
\smallskip

{In particular we  take $t=T$.
The estimate \eqref{eq:ee} gives that $\tau_N$ increases to $T$ a.s. as $N\to \infty$.
Passing to the limit, as
$N\to \infty$,  (\ref{apriori-est}) holds.}
\smallskip

{Taking the limit, as $N\to \infty$, in equality
\eqref{energy:eq}
we deduce
\eqref{energy:eqlim}.
}

\smallskip

This ends the proof.
$\hfill \;\blacksquare $

\smallskip

\begin{cor}
\label{Lp-energy:est}
Assume hypotheses of Proposition \ref{energy:est} and $u_0\in L^p(\Omega;H)$. Then for any $p\geq 4$
\begin{align}
\label{Lp-apriori-est}
&\sup_{0<\nu\leq \nu_0}\sup_{n}\left\{
\mathbb{E}\left(\sup_{0\leq r \leq T}\|u^\nu_n(r)\|_{L^2}^p\right)+\nu\int_0^T\mathbb{E}\left(\|u^\nu_n(s)\|_{L^2}^{p-2}\|u^\nu_n(s)\|_V^2\right)\,ds\right\}\notag\\
&\leq C(p,f,\mathcal{Q},\nu_0)\left(\mathbb{E}\left(\|u_0\|_{L^2}^p\right)+1\right).
\end{align}
\end{cor}

\textbf{Proof}.
{ For each $N\in \mathbb{N}$, let us consider the stopping time
$\tau_N=\inf\{t\geq 0:\|u_n^\nu(t)\|_{L^2}\geq N\}\wedge T$}. Applying It\^o's formula to expression (\ref{energy:eq}) and function $g(z)=z^{p/2}$,
\begin{align}
\label{Lp-energy:eq}
&\|u^\nu_n({t \wedge \tau_N})\|_{L^2}^p+p\nu\int_0^{t \wedge \tau_N}\|u^\nu_n(s)\|_{L^2}^{p-2}\|\nabla u^\nu_n(s)\|_{L^2}^2\,ds\notag\\
&=\|u^\nu_n(0)\|_{L^2}^p+p\nu\int_0^{t \wedge \tau_N}\left(\|u^\nu_n(s)\|_{L^2}^{p-2}\int_{\Gamma}(\kappa-\alpha)u^\nu_n(s)\cdot u^\nu_n(s)\,d\mathcal{S}\right)ds\notag\\
&\quad+
p\int_0^{t \wedge \tau_N}\|u^\nu_n(s)\|_{L^2}^{p-2}\langle  f(s),u^\nu_n(s)\rangle \,ds+p\int_0^{t \wedge \tau_N}\|u^\nu_n(s)\|_{L^2}^{p-2}\langle \sqrt{\mathcal{Q}}\,dW(s),u^\nu_n(s)\rangle \notag\\
&\quad +
\int_0^{t \wedge \tau_N}\biggl\{\frac{p}{2}\,{\rm tr}(\mathcal{Q})\|u^\nu_n(s)\|_{L^2}^{p-2}+\frac{p}{2}(p-2)\|u^\nu_n(s)\|_{L^2}^{p-4}
\langle \mathcal{Q}u^\nu_n(s),u^\nu_n(s)\rangle\biggr\}ds.
\end{align}
Using (\ref{NS:v1}),
\begin{align}
\label{Lp:NS:v1}
&\|u^\nu_n(s)\|_{L^2}^{p-2}\int_{\Gamma}|\kappa-\alpha|u^\nu_n(s)\cdot u^\nu_n(s)\,d\mathcal{S}\notag\\
&\leq \varepsilon_1 \|u^\nu_n(s)\|_{L^2}^{p-2}\|\nabla u^\nu_n(s)\|_{L^2}^2+C(\varepsilon_1)\|u^\nu_n(s)\|_{L^2}^{p}.
\end{align}
\smallskip
By Cauchy-Schwarz's inequality, we also get
\begin{align}
\label{Lp:NS:v2}
&\|u^\nu_n(s)\|_{L^2}^{p-2}|\langle f(s),u^\nu_n(s)\rangle|\leq \|u^\nu_n(s)\|_{L^2}^{p-2}\|f(s)\|_{L^2}\|u^\nu_n(s)\|_{L^2}\notag\\
&\leq\|u^\nu_n(s)\|_{L^2}^{p-1}\|f(s)\|_{L^2}\leq \|f(s)\|_{L^2}\left(1+\|u^\nu_n(s)\|_{L^2}^p\right).
\end{align}

\smallskip

Applying first Burkholder-Davis-Gundy's inequality and next Young's inequality with $p'=q'=1/2$,
\begin{align}
\label{Lp:NS:v3}
&\mathbb{E}\left|\sup_{0\leq r\leq t} \left\{p\int_0^{t \wedge \tau_N}\|u^\nu_n(s)\|_{L^2}^{p-2}\langle \sqrt{\mathcal{Q}}\,dW(s),u^\nu_n(s)\rangle \right\}\right|\notag\\
&\leq pC_1\,\mathbb{E}\left(\int_0^{t \wedge \tau_N}{\rm tr}(\mathcal{Q})\|u^\nu_n(s)\|_{L^2}^{2p-2}\,ds\right)^{1/2}\notag\\
&=\mathbb{E}\left(2\varepsilon_2\sup_{0\leq r\leq t}\left\{\|u^\nu_n(r \wedge \tau_N)\|_{L^2}^p\right\}\frac{1}{2\varepsilon_2}\int_0^{t \wedge \tau_N} p^2C_1^2{\rm tr}(\mathcal{Q})\|u^\nu_n(s)\|_{L^2}^{p-2}\,ds\right)^{1/2}\notag\\
&\leq \varepsilon_2 \,\mathbb{E}\left(\sup_{0\leq r \leq t}\|u^\nu_n(r \wedge \tau_N)\|_{L^2}^{p}\right)+C(\varepsilon_2,p) \mathbb{E}\left(\int_0^{t \wedge \tau_N}{\rm tr}(\mathcal{Q})\|u^\nu_n(r)\|_{L^2}^{p-2}\,ds\right),
\end{align}
where $C(\varepsilon_2,p)\doteq \frac{p^2C_1^2}{4\varepsilon_2}$.

\smallskip
The last term of the right hand side of \eqref{Lp-energy:eq}
can be estimated by
$$
\frac{1}{2}p(p-1)\int_0^{t \wedge \tau_N} {\rm tr}(\mathcal{Q})\|u^\nu_n(s)\|_{L^2}^{p-2}\,ds.
$$
Finally, notice that we can estimate
\begin{align}
\label{Lp:NS:v4}
\int_0^{t \wedge \tau_N} {\rm tr}(\mathcal{Q})\|u^\nu_n(s)\|_{L^2}^{p-2}\,ds\leq \int_0^{t \wedge \tau_N} {\rm tr}(\mathcal{Q})\left(1+\|u^\nu_n(s)\|_{L^2}^p\right)\,ds.
\end{align}
\smallskip
Thus, using expression (\ref{Lp-energy:eq}) and estimates (\ref{Lp:NS:v1}), (\ref{Lp:NS:v2}), (\ref{Lp:NS:v3}) and (\ref{Lp:NS:v4}), following the arguments used in the proof of Proposition \ref{energy:est}, and taking $\varepsilon_1=1-(2p)^{-1}$ and $\varepsilon_2=1/2$, one can complete the proof of this Corollary.
$\hfill \;\blacksquare $

\smallskip

The next lemma gives an important monotonicity property of operator $F$ in order to prove the existence and uniqueness for the weak solution, according to the Definition \ref{definicao}, to Equation (\ref{NS:eqevol}). As we shall see, from the stochastic point of view it will be a strong solution. Concerning weak solutions for stochastic Navier-Stokes equations, in the stochastic sense,  we refer \cite{AC90} and the more recent paper \cite{Yo13} (see also the references therein).
\begin{lem}
\label{monotonicity:lem}
For a given $r>0$ we consider the following (closed) $L^4$-ball $B_r$ in the space $V$:
$$
B_r\doteq \left\{v\in V:\,\|v\|_{[L^4(\mathcal{O})]^2}\leq r\right\}.
$$
Then the nonlinear operator $u\mapsto F(t,u)$, $t\in [0,T]$, is monotone in the convex ball $B_r$, that is, for any $u\in V$, $v\in B_r$, there exists a positive constant $C\doteq C(\nu_0,\mathcal{O},\alpha)$, depending on $\nu_0$, the domain $\mathcal{O}$ and $\alpha$ such that
\begin{equation}
\label{monotonicity:est}
\left\langle F(t,u)-F(t,v),u-v\right \rangle \leq C\left(1+\frac{r^4}{\nu^3}\right)\|u-v\|_{L^2}^2.
\end{equation}
\end{lem}
\textbf{Proof}.
Taking into account the definition of the operator $\mathcal{A}$  \eqref{A:def}, we have
\begin{align*}
\left\langle F(u)-F(v),u-v\right \rangle
&+\nu \int_{\mathcal{O}}|\nabla(u-v)|^2\,dx\\
&=
-\left\langle \mathcal{B}(u)-\mathcal{B}(v),u-v\right \rangle+
\nu\int_{\Gamma}(k-\alpha)|u-v|^2\,d\mathcal{S}.
\end{align*}
As in \eqref{NS:v1}, we derive the inequality
\begin{align*}
\nu \int_{\Gamma}(k-\alpha)|u-v|^2\,d\mathcal{S}\leq
\frac{\nu}{2}\left\|\nabla(u-v)\right\|^2_{L^2}+C \nu
\|u-v\|^2_{L^2}
\end{align*}
where $C$ is a constant which depends of $\mathcal{O}$  and $\alpha$.

\smallskip

For more details see the proof of Lemma 2.4 in \cite{MS02} and  Proposition 2.2 in \cite{SS06}.
$\hfill \;\blacksquare $

\smallskip

Now we shall prove the path-wise uniquenes of Equation (\ref{NS:eqevol}).

\smallskip

\begin{prop}
\label{uniqueness}
Assume the hypotheses of Proposition \ref{energy:est}. Let $u^\nu$ be a solution of Equation (\ref{NS:eqevol}), that is, an adapted stochastic process $u^\nu(t,x,\omega)$ satisfying (\ref{NS:eqevol}) and such that
$$
u^\nu\in L^2\left(\Omega;\mathcal{C}(0,T;H)\cap L^2(0,T;V)\right)\cap [L^4\left(]0,T[\times \mathcal{O}\times \Omega\right)]^2.
$$
If $v^\nu$ is another solution of Equation (\ref{NS:eqevol}) as an adapted stochastic process in the space $\mathcal{C}(0,T;H)\cap L^2(0,T;V)$, then
$$
\|u^\nu(t)-v^\nu(t)\|^2_{L^2}\exp\left\{-2C\int_0^t\left(1+\frac{1}{\nu^3}\|u^\nu(s)\|^4_{[L^4(\mathcal{O})]^2}\right)\,ds\right\}\leq \|u^\nu(0)-v^\nu(0)\|^2_{L^2},
$$
with probability $1$, for any $0\leq t\leq T$,  where $C$ is the positive constant that appears in Lemma \ref{monotonicity:lem}. In particular $u^\nu=v^\nu$, if $v^\nu$ satisfies the same initial condition as $u^\nu$.
\end{prop}
\textbf{Proof}.
Using Lemma \ref{monotonicity:lem}, it follows the same arguments as those in the proof of Proposition 3.2 in \cite{MS02}.
We should mention that this idea to prove  the path-wise uniqueness for the two dimensional stochastic Navier-Stokes equation already appear in \cite{Sc97}.
$\hfill \;\blacksquare $

\smallskip

The existence of solution to Equation (\ref{NS:eqevol}) is given in the following proposition.
\begin{prop}
\label{existence}
Suppose the hypotheses of Corollary \ref{Lp-energy:est}. Then there exists an adapted process $u^\nu(t,x,\omega)$ such that
$$
u^\nu\in L^p\left(\Omega;\mathcal{C}(0,T;H)\right)\cap L^2\left(\Omega;L^2(0,T;V)\right)\cap [L^4\left(]0,T[\times \mathcal{O}\times \Omega\right)]^2,
$$
and verifying Equation (\ref{NS:eqevol}). Furthermore,
\begin{align}
\label{Lp-apriori-est:sol}
&\sup_{0<\nu\leq \nu_0}\mathbb{E}\left\{
\sup_{0\leq r \leq T}\|u^\nu(r)\|_{L^2}^p+
\nu\int_0^T\|u^\nu(s)\|_V^2\,ds+\nu\int_0^T\|u^\nu(s)\|_{L^2}^{p-2}\|u^\nu(s)\|_V^2\,ds\right\}\notag\\
&\leq C(p,f,\mathcal{Q},\nu_0)\left(\mathbb{E}\left(\|u_0\|_{L^2}^p\right)+1\right).
\end{align}
\end{prop}
\textbf{Proof}.
Borrowing the arguments of the proof of Proposition 3.3 in \cite{MS02} and using the a priori estimates (\ref{apriori-est}) and (\ref{Lp-apriori-est}) and Lemma \ref{monotonicity:lem}, the proof of this Proposition can be completed.
$\hfill \;\blacksquare $

\smallskip

In the following section we shall consider the vorticity equation associated with Equation (\ref{NS:eq}) in order to improve the estimates (\ref{Lp-apriori-est:sol}).
 More precisely, we shall estimate  the $L^p$- norms of the vorticity process $\xi^\nu$  by the initial data, independently of the viscosity.

\smallskip

\section{$L^p$ a priori estimates for the vorticity independent of $\nu$}
\label{sec4}


Set $\xi^\nu={\rm curl}\,u^\nu$. We apply the operator ${\rm curl}$ to Equation (\ref{NS:eq}), obtaining the following vorticity equation:
\begin{equation}
\label{NS:curleq}
\left\{\renewcommand{\arraystretch}{1.8}
\begin{array}{ll}
\frac{\partial \xi^\nu(t)}{\partial t}-\nu \Delta \xi^\nu(t)+(u^\nu(t) \cdot\nabla )\xi^\nu(t)= {\rm curl}\,f(t)+{\rm curl}(\sqrt{\mathcal{Q}}\,\,\dot{W}(t)) &\text{ in }]0,T[\times \mathcal{O},    \\
\xi^\nu(0)={\rm curl}\,u_0 &\text{ in } \mathcal{O},\\
\xi^\nu=(2\kappa-\alpha)u^\nu\cdot \mathbf{t} &\text{ on }]0,T[\times \Gamma
\end{array}
\right.
\end{equation}
Notice that
$$
{\rm curl}(\sqrt{\mathcal{Q}}\,dW)=\sum_{k=1}^\infty \lambda_k^{-m}{\rm curl}\,v_k\,d\beta_k.
$$

\smallskip

In the following we shall denote by
 $\tilde{H}$ the space $L^2(\mathcal{O})$ endowed with the $L^2-$norm. We use the same notation for the $L^2-$norm of vector functions and scalar functions.
\smallskip

In the space $\tilde{H}$ consider the operator $\tilde{\mathcal{A}}:\,D(\tilde{\mathcal{A}})\subset \tilde{H}\rightarrow \tilde{H}$ with domain
$D(\tilde{\mathcal{A}})=\{\zeta\in
L^2(\mathcal{O}):\;\Delta\zeta\in
L^2(\mathcal{O})
\}$,
defined by $\tilde{\mathcal{A}} \zeta=-\Delta \zeta$ for all $\zeta\in D(\tilde{\mathcal{A}})$.

\smallskip

Set
$$
\zeta_k=\frac{{\rm curl}\, v_k}{\|{\rm curl}\, v_k\|_{L^2}}.
$$
We recall that the basis $\{v_k\}$ fixed previously was constructed in  \cite{CMR98}
 verifying the properties that
$\{{\rm curl}\,v_k\}$ is orthogonal in $L^2(\mathcal{O})$ and  for each $k$,
${\rm curl}\,v_k\in\mathcal{W}$ is an eigenfunction of the operator $\tilde{\mathcal{A}}$ with eigenvalue $\lambda_k.$ Then
the sequence
 $\{\zeta_k\}$ is an orthonormal  basis for the space $\tilde{H}$,  that verifies $\tilde{\mathcal{A}} \zeta_k=\lambda_k \zeta_k$. Thus,
$$
{\rm curl}(\sqrt{\mathcal{Q}}\,dW)=\sum_{k=1}^\infty \lambda_k^{-m}{\rm curl}\,v_k\,d\beta_k=\sum_{k=1}^\infty \lambda_k^{-m}\|{\rm curl}\,v_k\|_{L^2}\zeta_k\,d\beta_k.
$$

\smallskip

We define $\tilde{\mathcal{Q}}\in L(\tilde{H},\tilde{H})$ by
$$
\tilde{\mathcal{Q}} \zeta_k=\lambda_k^{-2m}\mu_k^2\zeta_k,
$$
where $\mu_k=\|{\rm curl}\,v_k\|_{L^2}$, and $\tilde W=\sum_{k=1}^\infty \zeta_k \beta_k$ is a new cylindrical Wiener process in $\tilde{H}.$

\smallskip

Notice that
\begin{equation}
\label{Ortho:prop1}
\|{\rm curl}\,v_k\|_{L^2}^2\leq C(1+\lambda_k)\|v_k\|_{L^2}^2.
\end{equation}

\smallskip

Indeed, (\ref{Ortho:prop1}) is a consequence of the following fact. We consider the following spectral problem that appears in the proof of Lemma 2.2 in \cite{CMR98}:
\begin{equation}
\label{spectral}
\left\{\renewcommand{\arraystretch}{1.8}
\begin{array}{ll}
\Delta^2 \psi=-\lambda\psi &\text{ in }\mathcal{O},  \\
-\Delta \psi=-(2\kappa-\alpha)\nabla\psi\cdot \mathbf{n} &\text{ on }\Gamma, \\
\psi=0 &\text{ on }\Gamma.
\end{array}
\right.
\end{equation}
Its variational form reads: find $\psi\in H^2(\mathcal{O})\cap H_0^1(\mathcal{O})$ and $\lambda\neq 0$ such that
$$
\int_\mathcal{O} \Delta \psi\, \Delta \varphi \,d x-\int_\Gamma (2\kappa-\alpha)\nabla \psi\cdot \mathbf{n}\,\nabla \varphi \cdot \mathbf{n}\,d\mathcal{S}=\lambda \int_\mathcal{O} \nabla \psi\, \nabla \varphi\,d x,\qquad \forall \varphi \in H^2(\mathcal{O})\cap H_0^1(\mathcal{O}).
$$
Finally, notice that
$$
v_k=-\nabla^\bot \psi_k:=(\partial_2 \psi_k,-\partial_1 \psi_k) \text{ and }{\rm curl}\,v_k=-\Delta \psi_k,
$$
for some $\psi_k$ solution of the spectral problem (\ref{spectral}).

\smallskip

Since $\sum_{k=1}^\infty \lambda_k^{-2m+1}<\infty$ we obtain that $\tilde{\mathcal{Q}}$ is a trace class operator.

\smallskip

Hence $\tilde{\mathcal{Q}}^{1/2}\tilde W$ is an $\tilde{H}$-valued centered Wiener process on $(\Omega,\mathcal{F},\mathbf{P})$, with covariance $\tilde{\mathcal{Q}}$ in $\tilde{H}$.

\smallskip

In terms of $\tilde{\mathcal{A}}$ and $\tilde{\mathcal{Q}}^{1/2}\tilde W$ we can write Equation (\ref{NS:curleq}) as
\begin{equation}
\label{NS:curleq2}
\left\{\renewcommand{\arraystretch}{1.8}
\begin{array}{ll}
d \xi^\nu(t)+\left\{\nu \tilde{\mathcal{A}} \xi^\nu(t)+(u^\nu(t) \cdot\nabla )\xi^\nu(t)\right\}dt\\
\hspace{5cm}= {\rm curl}\,f(t)\,dt+\tilde{\mathcal{Q}}^{1/2}\,d\tilde{W}(t) &\text{ in }]0,T[\times \mathcal{O},    \\
\xi^\nu(0)={\rm curl}\,u_0 &\text{ in } \mathcal{O},\\
\xi^\nu=(2\kappa-\alpha)u^\nu\cdot \mathbf{t} &\text{ on }]0,T[\times \Gamma
\end{array}
\right.
\end{equation}

\smallskip
The following  Lemma establishes a useful estimate for the elements  of the basis $\{v_j\}.$

\begin{lem}
\label{lem:basis}
Let $\{v_j\}$ be the previous fixed basis for $V$. Set $\xi_j={\rm curl}\,v_j$. Then
$$
\|\xi_j\|_{H^1(\mathcal{O})}\leq C(\lambda_j+1 )\|\xi_j\|_{L^2(\mathcal{O})}.
$$
\end{lem}
\textbf{Proof}.
We know that $\xi_j$ is solution of the Dirichlet problem
\begin{equation*}
\left\{\renewcommand{\arraystretch}{1.8}
\begin{array}{ll}
-\Delta \xi_j=\lambda_j\,\xi_j
\quad
 \text{ in }\Omega , \\
\xi_j=(2k-\alpha)v_j\cdot\mathbf{t} ,   \quad \text{ on }\Gamma
\end{array}
\right.
\end{equation*}
The functions $\xi_j$ can be written in the form
 $\xi_j=h_j+g_j$, where $h_j$ and $g_j$ verify
 \begin{equation*}
\left\{\renewcommand{\arraystretch}{1.8}
\begin{array}{ll}
-\Delta h_j=\lambda_j\,\xi_j  & \quad \text{ in }\mathcal{O} , \\
h_j=0 & \quad \text{ on }\Gamma
\end{array}
\right. \qquad \text{and}\qquad \left\{\renewcommand{\arraystretch}{1.8}
\begin{array}{ll}
-\Delta g_j=0 & \quad \text{ in }\mathcal{O} , \\
g_j=(2k-\alpha)v_j\cdot\mathbf{t}& \quad \text{ on }\Gamma
\end{array}
\right.
\end{equation*}
The functions $h_j$  and $g_j$ satisfy the
Calderon-Zygmund's estimates (see for example Theorems 1.8, 1.10 on pages 12, 15 and Proposition 1.2,
p. 14 in Girault and  Raviart \cite{GR86}))
\begin{equation*}
\|h_j\|_{H^{2}(\mathcal{O} )}\leq C\|\lambda_j\,\xi_j\|_{L^2(\mathcal{O} )},\quad \quad \|g_j\|_{H^1(\mathcal{O})}\leq C\|(2\kappa-\alpha)v_j\cdot\mathbf{t}\|_{H^{1/2}(\Gamma )}.
\end{equation*}
From trace's theory $\|v_j\|_{H^{1/2}(\Gamma )}\leq
\|v_j\|_{H^1(\mathcal{O} )}$.
On the other hand, we also know that $\psi_j$ verify
\begin{equation*}
\left\{\renewcommand{\arraystretch}{1.8}
\begin{array}{ll}
-\Delta {\psi}_j=\xi_j   \quad \text{ in }\mathcal{O} , \\
{\psi}_j=0  \quad \text{ on }\Gamma
\end{array}
\right.
\end{equation*}
and $v_j=-\nabla^{\perp} {\psi}_j$. Therefore
$$\|v_j\|_{H^1(\mathcal{O} )}\leq \|\psi_j\|_{H^2(\Omega)}\leq C\|\xi_j\|_{L^2(\mathcal{O})}.$$
Then we have
\begin{equation*}
\|\xi_j\|_{H^1(\mathcal{O})}\leq C\big(\lambda_j+1\big)\|\xi_j\|_{L^2(\mathcal{O} )}.
\end{equation*}
$\hfill \;\blacksquare $

The improvement on the a priori estimates obtained in Proposition \ref{energy:est} and
Corollary \ref{Lp-energy:est}  is given in the following result:

\begin{prop}
\label{existence-uniqueness:vorticity}
Suppose hypotheses of Proposition \ref{energy:est}. Assume also that $p>2$, ${\rm curl}\,f\in L^1\left(0,T;L^p(\mathcal{O})\right)$ and ${\rm curl}\,u_0\in L^p\left(\Omega;L^p(\mathcal{O})\right)$. Let $\xi^\nu$ be the vorticity
of $u^\nu$, then we have
\begin{align}
\label{better-Lp-apriori-est:sol}
&\sup_{\nu}\mathbb{E}\left(\sup_{0\leq r \leq T}\|\xi^\nu(r)\|_{L^p}^p\right)\notag\\
&\quad\leq  C\left({\rm curl}\,f,\tilde{\mathcal{Q}},T,p,\mathcal{O},\alpha\right)\left\{\mathbb{E}\left(\|u_0\|_{L^2}^p\right)+\mathbb{E}\left(\|{\rm curl}\,u_0\|_{L^p}^p\right)+1\right\}.
\end{align}
\end{prop}
\textbf{Proof}.
Let $u^\nu$ be a stochastic process which is solution of the stochastic Navier-Stokes equation \eqref{NS:eqevol} with
vorticity process $\xi^\nu$   solution of \eqref{NS:curleq2}.

Let us denote by $w$ the solution of the following linear  equation
\begin{equation}
\label{NS:curleq22}
\left\{\renewcommand{\arraystretch}{1.8}
\begin{array}{ll}
d w(t)+\left\{\nu \tilde{\mathcal{A}} w
(t)+(u^\nu(t) \cdot\nabla )w(t)\right\}dt=0 &\text{ in }]0,T[\times \mathcal{O},    \\
w(0)=0 &\text{ in } \mathcal{O},\\
w=(2\kappa-\alpha)u^\nu\cdot \mathbf{t} &\text{ on }]0,T[\times \Gamma
\end{array}
\right.
\end{equation}

We introduce the process $\rho=\xi^\nu-w.$
We can verify that $\rho$ is solution of the following stochastic differential equation:

 \begin{equation}
\label{NS:curleq21}
\left\{\renewcommand{\arraystretch}{1.8}
\begin{array}{ll}
d \rho(t)+\left\{\nu \tilde{\mathcal{A}} \rho(t)+(u^\nu(t) \cdot\nabla )\rho(t)\right\}dt\\
\hspace{5cm}= {\rm curl}\,f(t)\,dt+\tilde{\mathcal{Q}}^{1/2}\,d\tilde{W}(t) &\text{ in }]0,T[\times \mathcal{O},    \\
\rho(0)={\rm curl}\,u_0 &\text{ in } \mathcal{O},\\
\rho=0 &\text{ on }]0,T[\times \Gamma
\end{array}
\right.
\end{equation}
\smallskip

Using minor adaptation of the proof of Lemma 3 in \cite{LFLP05}, for $p>2$ we obtain that the  solution to Equation (\ref{NS:curleq22}) satisfies

\begin{align}
\label{curl:est1}
\|w\|_{L^\infty(0,T);L^p(\mathcal{O}))}\leq C(p,\mathcal{O},\alpha, \epsilon)\;\|u^\nu\|_{L^\infty(0,T);L^2(\mathcal{O}))}+\epsilon \|\xi^\nu\|
_{L^\infty(0,T);L^p(\mathcal{O}))},\quad \mathbb{P}\text{ -a.s.}-\omega
\end{align}
where $\epsilon$ is an arbitrary  small parameter.
Using Proposition \ref{energy:est}, we have
\begin{align}
\label{curl:est11}
\|w\|_{L^\infty(0,T);L^p(\mathcal{O}))}\leq C(p,\mathcal{O},\alpha, \epsilon)\;
\|u_0\|_{L^2(\mathcal{O})}+\epsilon \|\xi^\nu\|
_{L^\infty(0,T);L^p(\mathcal{O}))},\quad \mathbb{P}\text{ -a.s.}-\omega.
\end{align}

\smallskip

As regards Equation (\ref{NS:curleq21}),
let us denote by $\tilde{\mathcal{H}}$ the Cameron-Martin space of the $\tilde{H}$- valued Wiener process $\tilde{\mathcal{Q}}^{1/2} \tilde W$.
We introduce the class $R(\tilde{\mathcal{H}},L^p)$ of the so-called radonifying operators (see Definition 4.2 in \cite{BP01}, \cite{BP00} and \cite{BV03}).
Let $\{b_k\}$ be  a sequence of mutually independent  $N(0,1)$-distributed random variables and $\{h_k\}$ an orthonormal basis for $\tilde{\mathcal{H}}.$ The norm of an operator $K$  in this class of operators is defined by
$\|K\|_{R(\tilde{\mathcal{H}},L^p)}=\mathbb{E}\left(\left\|\sum_{k=1}^\infty b_k Kh_k\right\|_{L^p}^2
\right) .$
We remark that the notion of a radonifying operator is a generalization of the notion of a Hilbert-Schmidt operator to the case where $L^p$ is not a Hilbert space.
 In the particular case  $p=2$,  $R(\tilde{\mathcal{H}},L^p)$ is the space of Hilbert-Schmidt operators.

\smallskip

 We can verify that the inclusion
$I:\tilde{\mathcal{H}}\rightarrow L^p(\mathcal{O})$ belongs to the class
$R(\tilde{\mathcal{H}},L^p)$. In fact it is enough to verify that $I:\tilde{\mathcal{H}}\rightarrow H^1\cap L^p(\mathcal{O})$ is an Hilbert Schmidt operator (see Remark 6.1 in \cite{BL06} and Theorem 2.3 in \cite{BV03}).
 Considering the orthonormal basis $h_k=\tilde{\mathcal{Q}}^{1/2}\zeta_k$, the
Sobolev Imbedding Theorem (see Theorem 4.1.2 page 85 in \cite{AF03}), Lemma \ref{lem:basis} and estimates (\ref{Ortho:prop1}) and (\ref{Q:cond}), we obtain
   \begin{align*}
  \|I\|_{R(\tilde{\mathcal{H}},L^p)}^2&=\mathbb{E}\left(\left\|\sum_{k=1}^\infty \beta_k\tilde{\mathcal{Q}}^{1/2} \,\zeta_k\right\|_{L^p}^2\right)
\leq
\mathbb{E}\left(\left\|\sum_{k=1}^\infty \beta_k\tilde{\mathcal{Q}}^{1/2} \,\zeta_k\right\|_{H^1}^2\right)\\
&\leq
\mathbb{E}\left( \sum_{k=1}^\infty \beta_k^2 \left \langle\tilde{\mathcal{Q}}^{1/2} \,\zeta_k,\tilde{\mathcal{Q}}^{1/2} \,\zeta_k\right\rangle_{H^1} \right)+\mathbb{E}\left( 2\sum_{j<k} \beta_j\beta_k\left \langle\tilde{\mathcal{Q}}^{1/2} \,\zeta_j,\tilde{\mathcal{Q}}^{1/2} \,\zeta_k\right\rangle_{H^1} \right)\\
&= \sum_{k=1}^\infty \mathbb{E}\left(\beta_k^2\right) \left \langle\tilde{\mathcal{Q}} \,\zeta_k,\zeta_k\right\rangle_{H^1} + 2\sum_{j<k}\mathbb{E}\left( \beta_j\beta_k\right) \left \langle\tilde{\mathcal{Q}}^{1/2} \,\zeta_j,\tilde{\mathcal{Q}}^{1/2} \,\zeta_k\right\rangle_{H^1}\\
&=\sum_{k=1}^\infty \lambda_k^{-2m}\mu_k^2 \|\zeta_k\|_{H^1}^2=\sum_{k=1}^\infty \lambda_k^{-2m}\|{\rm curl}\,v_k\|_{H^1}^2\leq \mathcal{M}.
\end{align*}

Henceforth, $\langle \cdot,\cdot \rangle$ denotes the duality between $L^p$ and $L^{p/(p-1)}$ for some $1<p<\infty$.
For each $N\in\mathbb{N}$, we set
$\tau_N=\inf\{t\geq 0:\,\|\rho\|_{L^2}\geq N\}\wedge T.$ Taking the function $\Phi:L^p\rightarrow \mathbb{R},$ $\Phi(x)=\|x\|^p_{L^p}$
    and applying the It\^{o}'s formula  to the processes $\Phi(\rho(t))$
(see Theorem 4.3 of \cite{BP01}), we have
\begin{align}
\label{NS:curl:ito1}
&\|\rho(t\wedge\tau_N)\|_{L^p}^p=\|\rho(0)\|_{L^p}^p-p\int_0^{t\wedge\tau_N}\left\langle \nu\tilde{\mathcal{A}}\rho(s),|\rho(s)|^{p-2}\rho(s)\right\rangle\,ds \notag\\
&-p\int_0^{t\wedge\tau_N}\left\langle u^\nu_n(s)\cdot \nabla \rho(s),|\rho(s)|^{p-2}\rho(s)\right\rangle\,ds+p\int_0^{t\wedge\tau_N}\left\langle {\rm curl}\,f(s),|\rho(s)|^{p-2}\rho(s)\right\rangle \,ds \notag\\
&+p\int_0^{t\wedge\tau_N}\left\langle \tilde{\mathcal{Q}}^{1/2}\,d\tilde{W}(t),|\rho(s)|^{p-2}\rho(s)\right\rangle+\frac{1}{2}\int_0^{t\wedge\tau_N} {\rm tr}_{I}\,\Phi^{''}(\rho(s))\,ds,
\end{align}
where
\begin{equation*}
 {\rm tr}_{I}\,\Phi^{''}(v)\leq p(p-1)\|v\|_{L^p}^{p-2}\left\|I \right\|^2_{R(\mathcal{H},L^p)}\leq p(p-1)\|v\|_{L^p}^{p-2}\mathcal{M}.
\end{equation*}

\smallskip

Hence
\begin{align}
\label{curl:est2}
\frac{1}{2}\int_0^{t\wedge\tau_N} {\rm tr}_{I}\,\Phi^{''}(\rho(s))\,ds
&\leq \frac{1}{2}p(p-1)
\mathcal{M}
\int_0^{t\wedge\tau_N}\|\rho(s)\|_{L^p}^{p-2}\,ds\notag\\
&\leq \frac{1}{2}p(p-1)
\mathcal{M}
\int_0^{t\wedge\tau_N}\left(1+\|\rho(s)\|_{L^p}^p\right)\,ds.
\end{align}

\smallskip

Applying that $\rho=0$ on $\Gamma$, integration by parts formula and the fact that
\begin{equation}
\label{derivative}
\nabla\left[\left|\rho(s)\right|^{p-2}\rho(s)\right]=(p-1)\left|\rho(s)\right|^{p-2}\nabla \rho(s),
\end{equation}
we have
\begin{align}
\label{prop1}
\left\langle u^\nu(s)\cdot \nabla \rho(s),|\rho(s)|^{p-2}\rho(s)\right\rangle=0,
\end{align}

\smallskip

On the other hand, we consider the following identities and estimates for the remainder of the terms in (\ref{NS:curl:ito1}).
\begin{align}
\label{prop2}
&-p\int_0^{t\wedge\tau_N}\left\langle \nu\tilde{\mathcal{A}}\rho(s),|\rho(s)|^{p-2}\rho(s)\right\rangle\,ds \notag\\& =-\nu p(p-1)\int_0^{t\wedge\tau_N}\left(\int_\mathcal{O} \left|\nabla \rho(s,x)\right|^2\left|\rho(s,x)\right|^{p-2}\,dx\right)ds.
\end{align}
Indeed, using integration by parts formula and (\ref{derivative}),
\begin{align*}
&\left\langle \tilde{\mathcal{A}}\rho(s),\left|\rho(s)\right|^{p-2}\rho(s)\right\rangle=
\langle -\Delta \rho(s),|\rho(s)|^{p-2}\rho(s)\rangle\\
&=\int_\mathcal{O} \nabla \rho(s,x)\cdot \nabla\left[|\rho(s,x)|^{p-2}\rho(s,x)\right]d x\notag\\
&=(p-1)\int_\mathcal{O} \left|\rho(s,x)\right|^{p-2}\nabla \rho(s,x)\cdot \nabla\rho(s,x)\,d x\notag\\
&=(p-1)\int_\mathcal{O} \left|\nabla \rho(s,x)\right|^2\left|\rho(s,x)\right|^{p-2}\,d x.
\end{align*}

\smallskip

For the stochastic term, using Burkholder-Davies-Gundy inequality (see (6.10) page 1890 in \cite{BM10} or Theorem 4.2 in \cite{BP01}, for instance) and Young's inequality with $p'=q'=1/2$, we obtain
\begin{align}
\label{prop3}
&\mathbb{E}\left|\sup_{0\leq r \leq t}\left\{p\int_0^{r\wedge\tau_N}\left\langle \tilde{\mathcal{Q}}^{1/2}\,d\tilde{W}(t),|\rho(s)|^{p-2}\rho(s)\right\rangle\right\}\right|\notag\\
&\leq \varepsilon\,\mathbb{E}\left(\sup_{0\leq r \leq t}\left\|\rho(r\wedge\tau_N)\right\|_{L^p}^p\right)+\frac{ C_1^2 p^2\mathcal{M} }{4\varepsilon}\,\int_0^{t\wedge\tau_N}\left\{1+\mathbb{E}\left(\sup_{0\leq r \leq s}\left\|\rho(r)\right\|_{L^p}^p\right)\right\}ds.
\end{align}
In fact,
\begin{align*}
&\mathbb{E}\left|\sup_{0\leq r \leq t}\left\{p\int_0^{r\wedge\tau_N}\left\langle \tilde{\mathcal{Q}}^{1/2}\,d\tilde{W}(t),|\rho(s)|^{p-2}\rho(s)\right\rangle\right\}\right|\notag\\
&\leq C_1 p \, \mathbb{E}\left(\int_0^{t\wedge\tau_N} \|I\|_{R(\tilde{\mathcal{H}},L^p)}^2 \left\|\rho(s)\right\|_{L^p}^{2(p-1)}\,ds\right)^{1/2}\notag\\
&\leq C_1 p\,\mathbb{E}\left(\int_0^{t\wedge\tau_N}  \mathcal{M} \left\|\rho(s)\right\|_{L^p}^{2(p-1)}\,ds\right)^{1/2}\notag\\
&\leq \mathbb{E}\left(\left\{2\varepsilon\sup_{0\leq r \leq t}\left\|\rho(r\wedge\tau_N)\right\|_{L^p}^p\right\}^{1/2}\left\{\frac{1}{2\varepsilon}\int_0^{t\wedge\tau_N} C_1^2 p^2 \mathcal{M} \left\|\rho(s)\right\|_{L^p}^{p-2}\,ds\right\}^{1/2}\right)\notag\\
&\leq \varepsilon\,\mathbb{E}\left(\sup_{0\leq r \leq t}\left\|\rho(r\wedge\tau_N)\right\|_{L^p}^p\right)+ \frac{ C_1^2 p^2\mathcal{M} }{4\varepsilon}\,\mathbb{E}\left(\int_0^{t\wedge\tau_N} \left\|\rho(s)\right\|_{L^p}^{p-2}\,ds\right)\notag\\
&\leq \varepsilon\,\mathbb{E}\left(\sup_{0\leq r \leq t}\left\|\rho(r\wedge\tau_N)\right\|_{L^p}^p\right)+\frac{ C_1^2 p^2\mathcal{M} }{4\varepsilon}\,\int_0^{t\wedge\tau_N}\left\{1+\mathbb{E}\left(\sup_{0\leq r \leq s}\left\|\rho(r)\right\|_{L^p}^p\right)\right\}ds.
\end{align*}

\smallskip

Finally for the term with ${\rm curl}\,f$, applying H\"older inequality for $p>2$ and $q=p/(p-1)$
\begin{align}
\label{prop4}
\langle {\rm curl}\,f(s),|\rho(s)|^{p-2}\rho(s)\rangle &\leq \|{\rm curl}\,f(s)\|_{L^p}\||\rho(s)|^{p-1}\|_{L^{q}}\notag\\
&= \|{\rm curl}\,f(s)\|_{L^p}\|\rho(s)\|_{L^p}^{p-1}\notag\\
&\leq  \|{\rm curl}\,f(s)\|_{L^p}\left(1+\|\rho(s)\|_{L^p}^p\right).
\end{align}

\smallskip

To sum up, applying to (\ref{NS:curl:ito1}) the estimates (\ref{curl:est2}), (\ref{prop1}), (\ref{prop2}), (\ref{prop3}) and (\ref{prop4}), we obtain
\begin{align}
\label{NS:curl:ito2}
&\mathbb{E}\left(\sup_{0\leq r \leq t}\|\rho(r\wedge\tau_N)\|_{L^p}^p\right)+\nu p(p-1)\mathbb{E}\left(\int_0^{t\wedge\tau_N}\left(\int_\mathcal{O} |\nabla \rho(s,x)|^2|\rho(s,x)|^{p-2}\,dx\right)ds\right)\notag\\
&\leq \mathbb{E}\left(\|{\rm curl}\,u_0\|_{L^p}^p\right)+\int_0^{t\wedge\tau_N} p\,\|{\rm curl}\,f(s)\|_{L^p}\,ds\notag\\
&+\int_0^{t\wedge\tau_N} p\,\|{\rm curl}\,f(s)\|_{L^p}\,\mathbb{E}\left(\sup_{0\leq r \leq s}\|\rho(r\wedge\tau_N)\|_{L^p}^p\right)ds\notag\\
&+\varepsilon\,\mathbb{E}\left(\sup_{0\leq r \leq t}\|\rho(r\wedge\tau_N)\|_{L^p}^p\right)+\frac{ C_1^2 p^2\mathcal{M} }{4\varepsilon}\,\int_0^{t\wedge\tau_N}\left\{1+\mathbb{E}\left(\sup_{0\leq r \leq s}\|\rho(r\wedge\tau_N)\|_{L^p}^p\right)\right\}ds\notag\\
&+\frac{1}{2}p(p-1) \mathcal{M}\int_0^{t\wedge\tau_N}\left\{1+\mathbb{E}\left(\sup_{0\leq r \leq s}\|\rho(r\wedge\tau_N)\|_{L^p}^p\right)\right\}ds.
\end{align}

Finally, using (\ref{NS:curl:ito2}) with $\varepsilon=1/2$ and applying Gronwall-Bellman inequality for
\begin{align*}
X(t)&\doteq \mathbb{E}\left(\sup_{0\leq r \leq t}\|\rho(r)\|_{L^p}^p\right),\\
Y(t)&\doteq\nu p(p-1)\mathbb{E}\left(\int_0^{t\wedge\tau_N}\left(\int_\mathcal{O} |\nabla \rho(s,x)|^2|\rho(s,x)|^{p-2}\,dx\right)ds\right),\\
Z(t)&\doteq \mathbb{E}\left(\|{\rm curl}\,u_0\|_{L^p}^p\right)+\int_0^{t\wedge\tau_N} p\,\|{\rm curl}\,f(s)\|_{L^p}\,ds+t\left\{C_1^2p^2+p(p-1)\right\} \mathcal{M}/2,\\
\lambda(t)&\doteq Z(t)-Y(t),\\
\varphi(t)&\doteq p\,\|{\rm curl}\,f(t)\|_{L^p}+\left\{C_1^2p^2+p(p-1)\right\}\mathcal{M}/2\geq 0 \text{ with }\mathcal{K}\doteq \int_0^T \varphi(s)ds,
\end{align*}
we obtain
\begin{align*}
X(t)+Y(t)\leq Z(T)(1+\mathcal{K}e^{\mathcal{K}})\leq C(T,p,{\rm curl}\,f,\mathcal{M})\left(\mathbb{E}\left(\|{\rm curl}\,u_0\|_{L^p}^p\right)+1\right),
\end{align*}
uniformly in $n$ and $\nu$.

\smallskip

In particular, for $t=T$,
\begin{align*}
&
\mathbb{E}\left(\sup_{0\leq r \leq T}\|\rho(r\wedge\tau_N)\|_{L^p}^p\right)+2\nu p(p-1)\mathbb{E}\left(\int_0^{T\wedge\tau_N}\left(\int_\mathcal{O} |\nabla \rho(s,x)|^2|\rho(s,x)|^{p-2}\,d x\right)ds\right)\\
&\leq  C(T,p,{\rm curl}\,f,\mathcal{M})\left(\mathbb{E}\left(\|{\rm curl}\,u_0\|_{L^p}^p\right)+1\right),
\end{align*}
which gives that $\tau_N $ increses to $T$ a.s. as $N\to \infty.$ Taking the limit as
$N\to \infty$, we deduce

\smallskip

\begin{align}
\label{curl:est3}
\sup_{\nu}\mathbb{E}\left(\sup_{0\leq r \leq T}\|\rho(r)\|_{L^p}^p\right)\leq  C\left(T,p,{\rm curl}\,f,\tilde{\mathcal{Q}}\right)\left(\mathbb{E}\left(\|{\rm curl}\,u_0\|_{L^p}^p\right)+1\right).
\end{align}

\smallskip

Hence, estimates (\ref{curl:est11}) and (\ref{curl:est3}) yield the following estimate for the vorticity:

\begin{align}
\label{curl:est33}
\mathbb{E}\left(\|\xi^\nu\|_{L^\infty(0,T;L^p)}^p\right)&\leq C(p)
\mathbb{E}\left(\|\rho\|_{L^\infty(0,T;L^p)}^p+\|w\|_{L^\infty(0,T;L^p)}^p\right)\notag\\&\leq
 C\left(T,p,{\rm curl}\,f,\tilde{\mathcal{Q}}, \epsilon\right)\mathbb{E}\left(\|{\rm curl}\,u_0\|_{L^p}^p+\|\,u_0\|_{L^2}^p+1\right)\notag\\
 &\quad+C(p) \epsilon\;
 \mathbb{E}\left(\|\xi^\nu\|_{L^\infty(0,T;L^p)}^p\right).
\end{align}

Taking $\epsilon$ small enough we obtain
 \eqref{better-Lp-apriori-est:sol}.
$\hfill \;\blacksquare $

Using  Propositions \ref{energy:est} and \ref{existence-uniqueness:vorticity}, we can deduce the following result:
\begin{prop}
\label{properties:sequence}
Assume the hypotheses of Proposition \ref{existence-uniqueness:vorticity}. Then
\begin{equation}
\mathbb{E}\left(\|u^\nu\|^p_{L^\infty\left(0,T;[W^{1,p}(\mathcal{O})]^2\right)}\right)\leq C, \label{velocity:est1}
\end{equation}
with a constant $C>0$ independent of viscosity.
\end{prop}
\textbf{Proof}.
Owing to Poincar\'e's inequality, Lemma 3.1 in \cite{KP86} and a priori estimates (\ref{better-Lp-apriori-est:sol}) for the vorticity of $u^\nu$, (\ref{velocity:est1}) holds.
$\hfill \;\blacksquare $
\smallskip

\section{Vanishing viscosity limit}
\label{sec5}
In this section we shall prove our main result (Theorem \ref{main-result}), that is, the sequence of solutions $\{u^\nu\}_{0<\nu\leq \nu_0}$ to Equation (\ref{NS:eq}) converges to a solution of the stochastic Euler equations with the same initial
velocity as viscosity vanishes.

To establish the existence of  solution for the stochastic Euler equations, we follow a path-wise approach similar to \cite{BF99}. In our problem,  the vorticity of the involved processes do not vanish at the boundary, so, we need to estimate the boundary terms  which increases the difficulty. To overcome  such difficulties, we proceed analogously to deterministic methods in the articles \cite{CMR98}, \cite{LFLP05}.
Since the estimates, independent of the viscosity, are based on the maximum principle, we need to consider a regular Wiener process. To be precise, in the next two Lemmas the  Wiener process
$\sqrt{\mathcal{Q}}\,W(t)$ has covariance   $\mathcal{Q}=\mathcal{A}^{-2m}$,
 $m>4$.
\begin{lem}
\label{ex:Euler}
Assume that for a.e. $\omega\in \Omega $,
$u_0\in L^p(\mathcal{O})$ and  ${\rm curl\,}f\in
L^1(0,T;L^\infty(\mathcal{O}))$.  Let $u^\nu$ be the weak solution of \eqref{NS:eqevol}, then we have
\begin{align}
\label{velocity:pathwise:est}
&\|u^\nu(\omega)\|_{L^\infty\left(0,T;W^{1,p}(\mathcal{O})\right)}\leq C(\omega),\quad
&
\end{align}
where $C(\omega)$ does not depend on the viscosity $\nu$, for a.e. $\omega$
 in $ \Omega$ but depends on $\omega.$

Moreover, if we assume for a.e. $\omega\in \Omega $, $u_0\in L^\infty(\mathcal{O})$, the  estimate \eqref{velocity:pathwise:est} holds for $p=\infty$.
\end{lem}
\textbf{Proof}.
Let us consider the  stochastic process $v^\nu(t,\omega)\doteq u^{\nu}(t,\omega)-\sqrt{\mathcal{Q}}\,W(t, \omega)$
which satisfies  a deterministic equation similar to (16) in \cite{BF99}. By handling such equation we deduce
$\sup_{0\leq t\leq T}\|v^\nu(\omega)\|_{L^2(\mathcal{O})}\leq C$. To simplify, we represent by  $C$ a constant independent of the viscosity.

 Taking into account the regularity of the process  $\sqrt{\mathcal{Q}}\,W(t, \omega)$, we obtain
$$
\sup_{0\leq t\leq T}\|u^\nu(\omega)\|_{L^2(\mathcal{O})}\leq C.
$$
Next set $z^\nu(t,\omega)\doteq{\rm curl\,}u^{\nu}
(t,\omega)-{\rm curl}\left(\sqrt{\mathcal{Q}}\,W(t, \omega)\right)$. We consider just the case  $2<p<\infty$, since for $p=\infty$ is easier, the result follow directly by the maximum principle.
 For a.e. $\omega \in \Omega$, the sample paths of the process $z^\nu(t)$  verify, in the sense of the distributions  , the following equation:
\begin{equation}
\label{pointwise:vorticity}
\left\{\renewcommand{\arraystretch}{1.8}
\begin{array}{ll}
\frac{\partial z^\nu(t)}{\partial t}-\nu \Delta z^\nu(t)+(u^\nu(t) \cdot\nabla )
z^\nu(t)=g(t) &\text{ on }]0,T[\times \mathcal{O},  \\
z^\nu(0)=z_0 &\text{ in } \mathcal{O},\\
z^\nu=(2\kappa-\alpha)v^\nu\cdot \mathbf{t} &\text{ on }]0,T[\times \Gamma
\end{array}
\right.
\end{equation}
where
$g(t)\doteq {\rm curl}\,f(t)-
(u^\nu(t) \cdot\nabla )
\mathcal{Z}(t)+\nu \Delta \mathcal{Z}(t)$
and
 $\mathcal{Z}(t)\doteq{\rm curl}\left(\sqrt{\mathcal{Q}}\,W(t)\right)$.
 Let us denote by $\lambda\doteq\|(2\kappa-\alpha)v^\nu\cdot \mathbf{t}\|_{L^\infty([0,T]\times\Gamma)}$ and $L(t)\doteq\|g(t)\|
 _{L^\infty(\mathcal{O})}$. Given $u^\nu(t)$ , the linear problem
\begin{equation}
\label{pointwise:vorticity1}
\left\{\renewcommand{\arraystretch}{1.8}
\begin{array}{ll}
\frac{\partial \bar{z}(t)}{\partial t}-\nu \Delta \bar{z}(t)+(u^\nu(t) \cdot\nabla )
\bar{z}(t)= L(t)&\text{ on }]0,T[\times \mathcal{O},  \\
\bar{z}(0)=|z_0| &\text{ in } \mathcal{O},\\
\bar{z}=\lambda&\text{ on }]0,T[\times \Gamma
\end{array}
\right.
\end{equation}
is well posed with solution $\bar{z}(t)\in L^2(0,T;H^1(\mathcal{O})).$
Since  the function
  $z(t)=z^\nu(t)- \bar{z}(t)$
verifies the inequality
$$
\frac{\partial z(t)}{\partial t}-\nu \Delta z(t)+(u^\nu(t) \cdot\nabla )
z(t)\leq 0\quad\text{ on }]0,T[\times \mathcal{O},
$$
and is non positive on $[0,T ]\times \Gamma $ and at
$t=0$, \ the maximum principle implies that $z$ is a non positive
function, i.e.  $z^{\nu }\leqslant \bar{z}.$ Analogously, we
show that $z=-z^\nu-\bar{z}\leq  0.$ \ So, we conclude
that
 \begin{equation}
 \label{maximum}
 |z^\nu(t)|\leq \bar{z}(t), \quad\text{a. e. in }[0,t)\times \mathcal{O}.
 \end{equation}
 The difference process  $\hat{z}(t)\doteq\bar{z}(t)-\lambda$ verifies the  following equation:
\begin{equation}
\label{pointwise:vorticity2}
\left\{\renewcommand{\arraystretch}{1.8}
\begin{array}{ll}
\frac{\partial \hat{z}(t)}{\partial t}-\nu \Delta \hat{z}(t)+(u^\nu(t) \cdot\nabla )
\hat{z}(t)= L(t) &\text{ on }]0,T[\times \mathcal{O},\\
\hat{z}(0)=|z_0| -\lambda&\text{ in } \mathcal{O},\\
\hat{z}=0&\text{ on }]0,T[\times \Gamma.
\end{array}
\right.
\end{equation}
Multiplying the first equation of (\ref{pointwise:vorticity2}%
) by $G\doteq p|\hat{z}|^{p-2}\hat{z}$ and integrating over $\mathcal{O} ,$ we obtain
\begin{equation}
\frac{d}{dt}||\hat{z}||_{L_{p}(\mathcal{O}  )}^{p}+{\nu }p(p-1)\int_{\mathcal{O}
}|\hat{z}|^{p-2}|\nabla \hat{z}|^{2}\,dx\leqslant |\int_{\mathcal{O} }L(t)\ G\,dx|.  \label{a1}
\end{equation}%
Having%
\begin{eqnarray*}
\left|L(t)\int_{\Omega }\,G\,d x\right| &\leqslant &C\left[
||{\rm curl}f||_{L^\infty(\mathcal{O} )}+
||u
^{{\nu }}||_{C(\bar{\mathcal{O}})}
|| \mathcal{Z}||_{L{\infty}(\mathcal{O} )}
+
\nu|| \Delta \mathcal{Z}||_{L^{\infty}(\mathcal{O} )}
\right] ||\hat{z}||_{L^{p}(\mathcal{O} )}^{p-1}, \\
\end{eqnarray*}%
we verify that $||\hat{z}||_{L^{p}(\Omega )}^{p}$ satisfies a Bihari's type
inequality, that  gives the following estimate
\begin{equation*}
\Vert \hat{z}(t )\Vert _{L^{p}(\mathcal{O})}\leqslant C(\Vert z(0)\Vert
_{L^{p}(\mathcal{O} )}+\int_{0}^{t}
||v^{{\nu }}(r)||_{C(\overline{\mathcal{O} })}\,dr+1).
\end{equation*}%
Considering Nirenberg-Gagliardo's
 interpolation inequality%
\begin{equation*}
||v^{{\nu }}(t)||_{L^{\infty }(\mathcal{O} )}\leqslant C\left( ||v^{{\nu }}(t)||_{L^{2}(\mathcal{O} )}^{1-\theta }\Vert
v^{{\nu }}(t)\Vert
_{W^{1,p}(\mathcal{O} )}^{\theta }+||v^{{\nu }}(t)||_{L^{2}(\mathcal{O}
)}\right) ,\;\quad \theta =\frac{p}{2(p-1)},  \label{nir}
\end{equation*}%
the embedding theorem%
\begin{equation*}
W^{1,p}(\mathcal{O} )\hookrightarrow C^{\alpha }(\overline{\mathcal{O} })\qquad
\text{with\quad }\alpha =1-2/p  \label{emb}
\end{equation*}%
 we can write%
\begin{equation*}
||v^{{\nu }}(t)||_{C(\overline{\mathcal{O} })}=||v^{{\nu }}(t)||_{L^{\infty }(\mathcal{O} )}\leqslant C(\Vert v^{{\nu }}(t)\Vert _{L^{2}(\mathcal{O} )}+\Vert z^{\nu }\Vert
_{L^{p}(\mathcal{O} )}).  \label{eq9}
\end{equation*}%
Combining with \eqref{maximum}  we derive the following Gronwall's inequality for $z^\nu(t )$:
\begin{equation*}
||z^\nu(t)|| _{L^{p}(\mathcal{O} )}\leqslant C(\Vert z(0)\Vert
_{L^{p}(\mathcal{O})}+\int_{0}^{t}
||z^\nu(r,\cdot )|| _{L^{p}(\mathcal{O} )}\,dr+1),
\end{equation*}%
which implies $||z^\nu||_{L^\infty(0,T;{L^{p}(\mathcal{O} ))}}\leq C$, where $C$ is a constant independent of the viscosity.
Therefore, we have
\begin{equation}
||{\rm curl}\,u^\nu||_{L^\infty(0,T;{L^{p}(\mathcal{O} ))}}\leq C
\end{equation}
and consequently \eqref{velocity:pathwise:est} holds.

$\hfill \;\blacksquare $

\begin{lem}
\label{exist:Euler}
Under the assumptions of Lema \ref{ex:Euler}. Then  exists a stochastic process $u$ with sample paths in $C([0, T]; H)\cap L^\infty(0,T;W^{1,p}(\mathcal{O}))$, $p>2$ that is solution of the Euler equation \eqref{E:eqevol}, in the sense of (\ref{E:eqweak:time}). Moreover, in the case
$p=\infty$, such solution is unique.
\end{lem}
\textbf{Proof}.

Using estimates (\ref{velocity:pathwise:est}) and borrowing the arguments of Theorem 1 in \cite{LFLP05}
and Theorem 1.1 of \cite{BF99}, it can be proven the existence of $ u$, which is a  solution of \eqref{E:eqevol} in the sense of (\ref{E:eqweak:time}). To obtain a measurable
solution $ u$, we can use a measurable selection theorem (see Chapter 5 in \cite{Sr98} and also Lemma 3.1 in \cite{BF99} for a more specific result).
The proof of the uniquesess is standard.

$\hfill \;\blacksquare $

\smallskip

Finally, we can already prove our main result:

\textbf{Proof of Theorem \ref{main-result}.}
Notice that (i) is a consequence of Propositions \ref{uniqueness} and \ref{existence}.

Regarding (ii), observe that a stochastic process being solution of the Euler equation already exists, from
Lemma \ref{exist:Euler}. Let us suppose  $u_0(\omega)\in L^\infty(\mathcal{O})$ for a.e. $\omega \in \Omega$ and consider   $u^\nu$ and $u$   the unique solutions to Navier-Stokes equations  and  Euler equations, respectively. It remains to prove that $u^\nu$ converges to $u$ as the viscosity goes to zero.

Let us consider the difference process $u^\nu-u.$
\begin{align*}
\left\langle \mathcal{B}(u^\nu)-\mathcal{B}(u),u^\nu-u\right \rangle
&=  \label{equivalence-UniqnessWF1}
\int_{\mathcal{O}}\left[(u^\nu\cdot\nabla)(u^\nu-u)+((u^\nu-u)%
\cdot\nabla)u\right] \cdot\,(u^\nu-u)\; d x\\
&=
\int_{\Gamma} (u\cdot \mathbf{n})\; \,\frac{|u^\nu-u|^2}{2}\,dS+ \int_{\mathcal{O}} \left((u^\nu-u)%
\cdot\nabla\right)u\cdot (u^\nu-u)\,d x.\\
\end{align*}
Since $u\cdot \mathbf{n}=0$, the first term in the right hand side is zero.
 For the  second term we have
\begin{align*}
\left|\int_{\mathcal{O}}\left((u^\nu-u)%
\cdot\nabla\right)u\cdot (u^\nu-v)\,dx\right|
 &\leq \left\|\nabla u\right\|_{L^\infty}\|u^\nu-u\|^2_{L^2}.
 \end{align*}
On the other hand, taking into account the definition of the operator
$\mathcal{A}$ in \eqref{A:def}, We have
\begin{align*}
\left\langle \mathcal{A}(u^\nu),u^\nu-u\right \rangle
= \nu \int_{\mathcal{O}}\nabla u^\nu\cdot\nabla(u^\nu-u)\,dx
- \nu \int_{\Gamma}(k-\alpha)u^\nu\cdot(u^\nu-u)\,dS.
\end{align*}
Therefore, the difference process verifies  the following Gronwall inequality:
\begin{align*}
\frac{\partial}{\partial t}\|u^\nu(t)-u(t)\|^2_{L^2}\leq  C\nu +\|\nabla u(t)\|_{L^\infty}\|u^\nu(t)-u(t)\|^2_{L^2}
\end{align*}
for a.e. $\omega\in \Omega$,
which implies
$$\|u^\nu(t)-v(t)\|^2_{L^2}\leq C\nu\;\text{e}^{\int_0^t\|\nabla u(s)\|_{L^\infty}ds} $$
where $C$ is a constant independent of $\nu$. Then
 $\displaystyle \sup_{t\in[0,T]}\|u^\nu(t)-u(t)\|^2_{L^2}\to 0$
, as $\nu\to 0.$

\smallskip

$\hfill \;\blacksquare $

\smallskip




\smallskip
{\noindent \bf Acknowledgments:}
The first author thanks the support from FCT through the projects PTDC/MAT/104173/2008 and PEst-OE/MAT/UI0208/2011.
The second author wishes to acknowledge the support of the FCT portuguese project Pest-OE/MAT/UI0208/2011 by means of a post-doctoral position for one year.

\smallskip

The authors are grateful to
the referees for the  relevant suggestions and remarks that improved the article.

\end{document}